\documentclass{amsart}
\usepackage[latin1]{inputenc}
\usepackage{amsmath}
\usepackage{amsthm}
\usepackage{amsfonts}
\usepackage{amssymb}
\newtheorem{theorem}{Theorem}[section]
\newtheorem{lemma}[theorem]{Lemma}
\newtheorem{proposition}[theorem]{Proposition}
\newtheorem{definition}[theorem]{Definition}

\newtheorem{conjecture}[theorem]{Conjecture}
\newtheorem{remark}[theorem]{Remark}
\newtheorem{corollary}[theorem]{Corollary}

\begin{document}
\title[Spherical completeness of generalized numbers]{Spherical completeness of the non-archimedean ring of Colombeau generalized numbers}
\author[E. Mayerhofer]{Eberhard Mayerhofer}
\address{University of Vienna, Faculty of Mathematics, Nordbergstrasse 15, 1090 Vienna, Austria}
\thanks{Work supported by FWF research grants P16742-N04 and Y237-N13}
\email{eberhard.mayerhofer@univie.ac.at} 
\maketitle
\begin{abstract}
We show spherical completeness of the ring of Colombeau generalized
real (or complex) numbers endowed with the sharp norm. As an
application, we establish a Hahn-Banach extension theorem for
ultra-pseudo-normed modules (over the ring of generalized numbers) of generalized functions in the sense of
Colombeau.
\end{abstract}
\section{Introduction}
Let $(M,d)$ be an ultrametric space. For given $x\in M,\, r\in \mathbb
R^+$, we call $B_{\leq r}(x):=\{y\in M\mid d(x,y)\leq r\}$ the
dressed ball with center $x$ and radius $r$. Throughout $\mathbb N:=
\{1,2,\dots\}$ denote the positive integers. Let $(x_i)_i\in
M^{\mathbb N}$ and $(r_i)_i$ be a sequence of positive reals. We
call $(B_i)_i, \;B_i:=B_{\leq r_i}(x_i)\;(i\geq 1)$ a nested
sequence of dressed balls, if $r_1\geq r_2\geq r_3\dots$ and
$B_1\supseteq B_2\supseteq\dots$ . Following standard ultrametric
literature (cf.\ \cite{zAR}), nested sequences of dressed balls might have
an empty intersection. The converse property is defined as follows:
\begin{definition}\label{defsph}
$(M,d)$ is called spherically complete, if every nested sequence of
dressed balls has a non-empty intersection.
\end{definition}
It is evident that any spherically complete ultrametric space is
complete with respect to the topology induced by its metric (using
the well known fact that topological completeness of $(M,d)$ is
equivalent to the property of Definition \ref{defsph} with radii
$r_i \searrow 0$) . However, there are popular non-trivial examples
in the literature, for which the converse is not true. As an example
we mention the completion $\mathbb C_p$ of the algebraic
closure of the field  of rational $p$-adic numbers. Due
to Krasner, this field has nice algebraic properties (as it is
algebraically closed, and even isomorphic to the complex numbers
cf.\ \cite{zAR}, pp.\ 134--145), but it also has been shown, that
$\mathbb C_p$ is not spherically complete. This is mainly due to the
fact that the complex $p$-adic numbers are a separable, complete
ultrametric space with dense valuation (cf.\ \cite{zAR}, pp.\
143--144). However, for an ultrametric field $K$, spherical
completeness is necessary in order to ensure $K$ has the Hahn-Banach
extension property (to which we refer as HBEP), that is, any ultra-normed
$K$-vector space $E$ admits continuous linear functionals previously
defined on a strict subspace $V$ of $E$ to be extended to the whole
space under conservation of their norm (this is due to W. Ingleton,
\cite{Ingleton}). Since spherical completeness fails, it is natural
to ask if the $p$-adic numbers could at least be spherically
completed, i.e., if there existed a spherically complete ultrametric
field $\Omega$ into which $\mathbb C_p$ can be embedded. This
question has a positive answer (cf.\ \cite{zAR}). The necessity of
spherical completeness for the HBEP of $K=\mathbb C_p$ is evident:
even the identity map
\[
\varphi:\;\; \mathbb C_p\rightarrow\mathbb C_p,\quad \varphi(x):=x
\]
cannot be extended to a functional $\psi:\Omega\rightarrow\mathbb
C_p$ under conservation of its norm $\|\varphi\|=1$
 (here we consider $\Omega$ as a $\mathbb C_p$- vector space).\footnote{To check this, let $B_i:=B_{\leq r_i}(x_i)$ be a nested sequence of dressed balls in $\mathbb C_p$ with empty intersection. Then $\hat B_i:=B_{\leq r_i}(x_i)\subseteq\Omega$ have nonempty intersection, say $\Omega\ni\alpha\in\bigcap_{i=1}^{\infty}\hat B_i$. Assume further, the identity $\varphi$ on $\mathbb C_p$ can be extended to some linear map $\psi:\Omega\rightarrow\mathbb C_p$ under conservation of its norm. Then
\[
\vert \psi(\alpha)-x_i\vert_{\Omega}=\vert
\psi(\alpha)-\varphi(x_i)\vert_{\Omega}\leq
\|\psi\|\vert\alpha-x_i\vert_{\mathbb
C_p}=\vert\alpha-x_i\vert_{\mathbb C_p},
\]
therefore $\psi(\alpha)\in \bigcap_{i=1}^{\infty}B_i$ which is a
contradiction and we are done.}

The present paper is motivated by the question if a HBEP for the ring $\widetilde{\mathbb R}$ (resp.\ $\widetilde{\mathbb C}$) of generalized numbers
holds. Even though a first version of Hahn-Banach's Theorem is given in (\cite{Garetto1}, Proposition 3.23), a general version of the latter has not been established yet in the literature.

The analogy with the $p$-adic case lies at hand, since the ring of generalized numbers can naturally 
be endowed with an ultrametric pseudo-norm.
However, the presence of zero-divisor in $\widetilde{\mathbb R}$ as
well as the failing multiplicativity of the pseudo-norm turns the
question into a non-trivial one and Ingleton's ultrametric version
of the Hahn-Banach Theorem cannot be carried over to our setting
unrestrictedly.

On our first step tackling this question we discuss spherical
completeness of the ring of generalized numbers endowed with the
given ultrametric (induced by the respective ultra-pseudo-norm, cf.\ the preliminary
section).

$\widetilde{\mathbb R}$ first was introduced as the set of values of
generalized functions at standard points; however, a subring
consisting of compactly supported generalized numbers turned out to
be the set of points for which evaluation determines uniqueness,
whereas standard points do not suffice do determine generalized
functions uniquely (cf. \cite{Eyb4, MO1} as well as section 1.2.4 in \cite{Bible}). A hint, that
$\widetilde{\mathbb R}$ (or $\widetilde{\mathbb C}$ as well), the
ring of generalized real (or complex) numbers is spherically
complete, is, that contrary to the above outlined situation on
$\mathbb C_p$, the generalized numbers endowed with the topology
induced by the sharp ultra-pseudo-norm are not separable. This, for
instance, follows from the fact that the restriction of the sharp
valuation (cf.\ Section \ref{presec}) to the real (or complex) numbers is discrete. 

Having motivated our work by now, we may formulate the aim of this paper,
which is to prove the following:
\begin{theorem}\label{theoremsph}
The ring of generalized numbers is spherically complete.
\end{theorem}
We therefore have an independent proof of the fact (cf.\
\cite{Garetto1}, Proposition 1.\ 31 and Proposition 3.4):
\begin{corollary}
The ring of generalized numbers is topologically complete.
\end{corollary}
In the last section of this paper we present a modified version of
Hahn-Banach's Theorem which bases on spherically completeness of
$\widetilde{\mathbb R}$ (resp.\ $\widetilde{\mathbb C}$). Finally, a remark on
the applicability of the ultrametric version of Banach's fixed point theorem can be found in the appendix.
\section{Preliminaries}\label{presec}
In what follows we repeat the definitions of the ring of (real or complex) generalized numbers along with its non-archimedean valuation function $\rm v$. The material is taken from different sources; as references we recommend the recent works due to C.\ Garetto (\cite{Garetto2,Garetto1}) and A.\ Delcroix et al (\cite{DHPV}) as well as one of the original sources of this topic due to D.\ Scarpalezos (cf.\ \cite{Scarpalezos1}).\\
Let $I:=(0,1]\subseteq \mathbb R$, and let $\mathbb K$ denote
$\mathbb R$ resp.\ $\mathbb C$. The ring of generalized numbers over
$\mathbb K$ is constructed in the following way: given the ring of
moderate (nets of) numbers
\[
\mathcal E_M:=\{(x_{\varepsilon})_{\varepsilon}\in\mathbb K^I\mid
\exists\; N:\vert
x_{\varepsilon}\vert=O(\varepsilon^{-N})\;(\varepsilon\rightarrow 0)\}
\]
and, similarly, the ideal of negligible nets in $\mathcal E_M$ which are of the form
\[
\mathcal N:=\{(x_{\varepsilon})_{\varepsilon}\in\mathbb K^I\mid
\forall\; m:\vert
x_{\varepsilon}\vert=O(\varepsilon^m)\;(\varepsilon\rightarrow 0)\},
\]
we define the generalized numbers as the factor ring
$\widetilde{\mathbb K}:=\mathcal E_M/\mathcal N$. We define a valuation function $\rm v$ on $\mathcal E_M$ 
with values in $(-\infty,\infty]$ in the following way:
\[
{\rm v}((u_{\varepsilon})_{\varepsilon}):=\sup\,\{b\in\mathbb R\mid
\vert
u_{\varepsilon}\vert=O(\varepsilon^b)\;\;(\varepsilon\rightarrow
0)\}.
\]
This valuation can be carried over to the ring of generalized
numbers in a well defined way, since for two representatives of a
generalized number, the valuation above coincides (cf.\
\cite{Garetto1}, Section 1). We then endow $\widetilde{\mathbb
K}$ with an ultra-pseudo-norm ('pseudo' refers to
non-multiplicativity) $\vert \;\; \vert_{\rm e}$ in the following way:
$\vert 0\vert_{\rm e}:=0$, and whenever $x\neq 0$, $\vert
x\vert_{\rm e}:={\rm e}^{-{\rm v}(x)}$. With the ultrametric $d_{\rm e}$ induced by the above
ultra-pseudo-norm, $\widetilde{\mathbb K}$ turns out to be a non-discrete
ultrametric space, with the following topological properties:
\begin{enumerate}
\item $(\widetilde{\mathbb K},d_{\rm e} )$ is topologically complete (cf. \cite{Garetto1}),
\item $(\widetilde{\mathbb K},d_{\rm e})$ is not separable, since the restriction of $d_{\rm e}$ onto $\mathbb K$ is discrete.
\end{enumerate}
The latter property holds, since on metric spaces second
countability and separability are equivalent and the well known fact
that the property of second countability is inherited by subspaces
(whereas separability is not in general).

In order to avoid confusion we henceforth denote closed balls in
$\mathbb K$ by $B_{\leq r}(x):=\{y\in\mathbb K\;|\;\vert y-x\vert\leq r\}$ in distinction with dressed balls in
$\widetilde{\mathbb K}$ which we denote by $\widetilde B_{\leq
r}(x):=\{y\in\widetilde{\mathbb K}\;|\;\vert y-x\vert_{\rm e}\leq r\}$. Similarly stripped balls and the sphere in the ring of
generalized numbers are denoted by $\widetilde B_{< r}(x)
:=\{y\in\widetilde{\mathbb K}\;|\;\vert y-x\vert_{\rm e}< r\}$
 resp.\ $\widetilde S_r(x):=\{y\in\widetilde{\mathbb K}\;|\;\vert y-x\vert_{\rm e}= r\}$. \section{Euclidean Models of sharp
neighborhoods} Throughout, a net of real numbers
$(C_\varepsilon)_\varepsilon$ is said to {\it increase monotonically
with} $\varepsilon\rightarrow 0$, if the following holds:
\[
\forall \eta,\eta'\in I:\;(\eta\leq\eta'\Rightarrow C_\eta\geq
C_{\eta'}).
\]
To begin with we formulate the following condition:\\
{\it Condition (E).}\\
A net $(C_\varepsilon)_\varepsilon$ is said to
satisfy condition (E), if it is
\begin{enumerate}
\item positive for each $\varepsilon$ and
\item monotonically increasing with $\varepsilon\rightarrow 0$, and finally, if
\item the sharp norm is $\vert(C_\varepsilon)_\varepsilon\vert_e=1$.
\end{enumerate}
Next, we introduce the notion of euclidean models for sharp
neighborhoods of generalized points:
\begin{definition}\rm
Let $x\in\widetilde{\mathbb K}$, $\rho\in\mathbb R,\;
r:=\exp(-\rho)$. Let further
$(C_{\varepsilon})_{\varepsilon}\in\mathbb R^I$ be a net of real
numbers satisfying condition (E) and let
$(x_{\varepsilon})_{\varepsilon}$ be a representative of $x$. Then
we call the net of closed balls
$(B_{\varepsilon})_{\varepsilon}\subseteq \mathbb K^I$ given by
\[
B_{\varepsilon}:=B_{\leq
C_{\varepsilon}\varepsilon^{\rho}}(x_{\varepsilon})
\]
for each $\varepsilon\in I$ an euclidean model for $\widetilde
B_{\leq r}(x)$.
\end{definition}
Note, that every dressed ball admits an euclidean model: let
$(x_\varepsilon)_\varepsilon$ be a representative of $x$ and define
$(C_\varepsilon)_\varepsilon$ by $C_{\varepsilon}:=1$ for each
$\varepsilon\in I$; then $B_{\leq
C_{\varepsilon}\varepsilon^{\rho}}(x_{\varepsilon})$ determines an euclidean model for $\widetilde B_{\leq r}(x)$ when $\rho=-\log(r)$.\\We
need to mention that whenever we write
$(B^{(1)}_{\varepsilon})_{\varepsilon}\subseteq
(B^{(2)}_{\varepsilon})_{\varepsilon}$, we mean the inclusion
relation $\subseteq$ holds component wise (that is for each
$\varepsilon\in I$), and we say
$(B^{(2)}_{\varepsilon})_{\varepsilon}$ contains
$(B^{(1)}_{\varepsilon})_{\varepsilon}$.\\ The following lemma is
basic; however, in order to get familiar with the concept of
euclidean neighborhoods, we include a detailed proof:
\begin{lemma}\label{capture}
For $x\in\widetilde{\mathbb K}$ and $r>0$ let
$(B_{\varepsilon})_{\varepsilon}$ be an euclidean model for
$\widetilde B_{\leq r}(x)$. Then,
\begin{enumerate}
\item \label{capture1} for any $y\in \widetilde B_{<r}(x)$ there exists a representative $(y_{\varepsilon})_{\varepsilon}$ such that $y_{\varepsilon}\in B_{\varepsilon}$ for all $\varepsilon \in I$.
\item \label{capture2} There exist $y\in \widetilde S_r(x)$ fulfilling the following property: for every representative $(y_\varepsilon)_\varepsilon$
 of $y$ there exists $\varepsilon_0\in I$ such that $y_{\varepsilon_0}\not\in B_{\varepsilon_0}$. However,
 for all $y\in\widetilde S_r(x)$ and for all representatives $(y_\varepsilon)_\varepsilon$ of $y$ there exists
an euclidean model $\hat B_\varepsilon:=\hat B_{\leq \hat C_\varepsilon\varepsilon^\rho}(x_\varepsilon)$ for
$\widetilde B_{\leq r}(x)$ containing $(B_\varepsilon)_\varepsilon$ such that $y_\varepsilon\in\hat B_\varepsilon$
and $d(\partial \hat B_{\varepsilon},y_{\varepsilon})\geq \frac{C_{\varepsilon}}{2}\varepsilon^{\rho}$ for all $\varepsilon \in I$.
\end{enumerate}
\end{lemma}
\begin{proof}
(\ref{capture1}): By definition of the sharp norm, $\vert
y-x\vert_e<r$ is equivalent to the situation, that for each
representative $(y_{\varepsilon})_{\varepsilon}$ of $y$ and for each
representative $(x_{\varepsilon})_{\varepsilon}$ of $x$, we have
\[
\sup\{b\in\mathbb R\mid \vert y_{\varepsilon}-
x_{\varepsilon}\vert=O(\varepsilon^b) (\varepsilon\rightarrow
0)\}>\rho,
\]
and this implies that there exists some $\rho'>\rho$ such that for
any representative $(y_{\varepsilon})_{\varepsilon}$ of $y$ and any
representative $(x_{\varepsilon})_{\varepsilon}$ of $x$ we have
\[
\vert
y_{\varepsilon}-x_{\varepsilon}\vert=o({\varepsilon}^{\rho'}),\quad
\varepsilon\rightarrow 0.
\]
This further implies that for any choice of representatives of $x$
resp.\ of $y$, there exists some $\eta\in I$ with
\begin{equation}\label{tinyradius}
\vert y_{\varepsilon}-x_{\varepsilon}\vert\leq \varepsilon^{\rho'}
\end{equation}
for each  $\varepsilon<\eta$. Since $C_{\varepsilon}>0$ for each $\varepsilon\in I$ and $C_\varepsilon$ is monotonically increasing with $\varepsilon\rightarrow0$, we have $\varepsilon^{\rho'}\leq C_{\varepsilon}\varepsilon^{\rho}$ for sufficiently small $\varepsilon$. Therefore, a suitable choice of $\eta$ and of $y_{\varepsilon}$ for $\varepsilon\geq \eta$ yields the
first claim (for instance, one can set
$y_{\varepsilon}:=x_{\varepsilon}$ whenever $\varepsilon\geq \eta$).\\
We go on by proving (\ref{capture2}): For the first part, set
\[
y_{\varepsilon}:=2C_{\varepsilon}\varepsilon^{\rho}+x_{\varepsilon}
\]
Let $y$ denote the class of $(y_{\varepsilon})_{\varepsilon}$. It is
evident, that $y\in \widetilde S_r(x)$. However,
$(y_{\varepsilon})\notin B_{\varepsilon}$ for each $\varepsilon\in
I$. Indeed,
\[
\forall\;\varepsilon\in I:\vert
y_{\varepsilon}-x_{\varepsilon}\vert= 2
C_{\varepsilon}\varepsilon^{\rho}>C_{\varepsilon}\varepsilon^{\rho},
\]
since $C_\varepsilon>0$ for each $\varepsilon$. We further show,
that the same holds for any representative $(\bar
y_{\varepsilon})_{\varepsilon}$ of $y$ for sufficiently small index
$\varepsilon$. Indeed, the difference of two representatives being
negligible implies that for any $N>0$ we have
\[
y_{\varepsilon}-\hat y_{\varepsilon}=o(\varepsilon^N)\;\;
(\varepsilon\rightarrow 0).
\]
Therefore, for $N>\rho$ and sufficiently small $\varepsilon$, we
have:
\[
\vert \hat y_{\varepsilon}-y_{\varepsilon}\vert\geq \vert \vert\hat
y_{\varepsilon}-y_{\varepsilon}\vert - \vert
y_{\varepsilon}-x_{\varepsilon}\vert\vert\geq
2C_{\varepsilon}\varepsilon^{\rho}-\varepsilon^N\geq
\frac{3}{2}C_{\varepsilon}\varepsilon^{\rho}>C_{\varepsilon}\varepsilon^{\rho}.
\]
Therefore we have shown the first part of (\ref{capture2}). Let us take an arbitrary
$y\in \widetilde S_r(x)$. We demonstrate how to blow up
$(B_{\varepsilon})_{\varepsilon}$ to catch some fixed representative
$(y_{\varepsilon})_{\varepsilon}$ of $y$. Since $\vert
y-x\vert=e^{-\rho}=r$, there is a net $C'_{\varepsilon}\geq 0$
($\vert (C'_{\varepsilon})_{\varepsilon}\vert_e=1$) such that
\[
\forall\varepsilon\in I:\;\vert
y_{\varepsilon}-x_{\varepsilon}\vert=C_{\varepsilon}'\varepsilon^{\rho}
\]
Set $C''_{\varepsilon}=\max_{\eta\geq\varepsilon}\{1,C'_{\eta}\}$.
This ensures that $(C''_{\varepsilon})$ is a monotonically increasing
with $\varepsilon\rightarrow 0$, above $1$ for each $\varepsilon\in
I$, and $\vert (C''_{\varepsilon})\vert_e=1$ is preserved. The same holds for
the net $C'''_\varepsilon:=C''_\varepsilon+C_\varepsilon$.
Define $B'_{\varepsilon}:=B_{\leq
C'''_{\varepsilon}\varepsilon^{\rho}}(x_{\varepsilon})$.
Then $(B'_{\varepsilon})_{\varepsilon}$ is a new model for
$\widetilde B_{\leq r}(x)$ containing the old model and
$(y_{\varepsilon})_{\varepsilon}$ as well, since the sum
$C'''_{\varepsilon}$ satisfies the required properties (of condition (E)), and
\[
\vert y_{\varepsilon}-x_{\varepsilon}\vert\leq
C''_{\varepsilon}\varepsilon^{\rho}\leq
C_{\varepsilon}'''\varepsilon^{\rho}.
\]
Setting $\hat C_{\varepsilon}:=2 C'''_\varepsilon$ we obtain a model $\hat B_\varepsilon:=B_{\leq \hat C_\varepsilon\varepsilon^\rho}(x_\varepsilon)$ for $\widetilde B_{\leq r}(x)$ with the further property that $\vert
y_{\varepsilon}-x_{\varepsilon}\vert\leq
\frac{C_{\varepsilon}'''}{2}\varepsilon^\rho$ for each $\varepsilon\in I$ which
finishes the proof of (\ref{capture2}).\end{proof} 
\begin{remark}\rm
The preceding lemma can be reformulated in the following way: For
all $y\in\widetilde B_{\leq r}(x)$ there exists an euclidean model $B_\varepsilon:=
B_{\leq C_{\varepsilon}\varepsilon^{\rho}}(x_{\varepsilon})$ and a representative $(y_\varepsilon)_\varepsilon$ of $y$
such that $y_\varepsilon\in B_\varepsilon$ and $d(\partial B_\varepsilon,y_\varepsilon)\geq \frac{C_\varepsilon}{2}\varepsilon^\rho$
for all $\varepsilon\in I$.
\end{remark}
Before going on by establishing the
crucial statement which will allow us to translate decreasing
sequences of closed balls in the given ultrametric space
$\widetilde{\mathbb K}$ to decreasing sequences of their
(appropriately chosen) euclidean models, we introduce a useful term:
\begin{definition}\rm
Suppose, we have a nested sequence $(\widetilde B_i)_{i=1}^{\infty}$
of closed balls with centers $x_i$  and radii $r_i$ in
$\widetilde{\mathbb K}$. Let $(B^{(i)}_{\varepsilon})_{\varepsilon}$ be an euclidean model for $\widetilde B_i$ ($i\in\mathbb N$). We say that this associated sequence of euclidean models is proper, if $\left((B^{(i)}_{\varepsilon})_{\varepsilon}\right)_{i=1}^{\infty}$
is nested as well, that is, if we have:
\[
(B^{(1)}_{\varepsilon})_{\varepsilon}\supseteq
(B^{(2)}_{\varepsilon})_{\varepsilon}\supseteq(B^{(3)}_{\varepsilon})_{\varepsilon}\supseteq
\dots\,.
\]
\end{definition}
\section{Proof of the main Theorem} In order to prove the main
statement, we proceed by establishing two important preliminary
statements. First, a remark on the notation adopted in the sequel: if
$(x_i)_i$, a sequence of points in the ring of generalized numbers,
is considered, then $(x_\varepsilon^{(i)})_\varepsilon$ denote
(certain) representatives of the $x_i$'s. Furthermore, for
subsequent choices of nets of real numbers
$(C^{(i)}_\varepsilon)_\varepsilon$, and positive radii $r_i$, we
denote by $\rho_i$ the negative logarithms of the $r_i$'s
($i=1,2,\dots,$) while the euclidean models for the balls $\widetilde
B_{\leq r_i}(x_i)$ with radii $r_\varepsilon
^{(i)}:=C_\varepsilon^{(i)}\varepsilon^{\rho_i}$ to be constructed are
denoted by
\[
B_\varepsilon^{(i)}:=B_{\leq
r_\varepsilon^{(i)}}(x_\varepsilon^{(i)}).
\]
We start with the fundamental proposition:
\begin{proposition}\label{lemmatique}
Let $x_1,\;x_2\in\widetilde{\mathbb K}$, and $r_1,\;r_2$ be positive
numbers such that $\widetilde B_{\leq r_1}(x_1)\supseteq \widetilde
B_{\leq r_2}(x_2)$. Let $(x^{(1)}_\varepsilon)_\varepsilon$ be a representative of
$x_1$. Then the following holds:
\begin{enumerate}
\item \label{lemmatique1} There exists a net $(C^{(1)}_\varepsilon)_\varepsilon$ satisfying condition (E) and a representative
$(x_\varepsilon^{(2)})_\varepsilon$ of $x_2$ such that 
\begin{equation}\label{eqcontained}
x_\varepsilon^{(2)}\in B_{\leq
\frac{C_\varepsilon^{(1)}\varepsilon^{\rho_1}}{2}}(x_\varepsilon^{(1)})
\end{equation}
for each $\varepsilon\in I$.
\item \label{lemmatique2} Furthermore, for each net $(C^{(2)}_\varepsilon)_\varepsilon$ satisfying condition (E)
there exists $\varepsilon_0^{(1)}\in I$ such that $B^{(2)}_{\varepsilon}\subseteq B^{(1)}_{\varepsilon}$
for all $\varepsilon\in (0,\varepsilon_0^{(1)})$.
\end{enumerate}
\end{proposition}
\begin{proof}
Proof of (\ref{lemmatique1}): We distinguish the following two cases:
\begin{itemize}
\item $x_2\in \widetilde S_{r_1}(x_1)$, that is $\vert x_2-x_1\vert_{\rm e}=r_1$. For a given representative $(x_{\varepsilon}^{(2)})_{\varepsilon}$ of $x_2$, define $\hat C_{\varepsilon}^{(1)}:=\vert x_{\varepsilon}^{(1)}-x_{\varepsilon}^{(2)}\vert$. Now, set $ C_{\varepsilon}^{(1)}:=2\max(\{\hat C_{\eta}^{(1)} \vert \eta>\varepsilon\},1)$. Then not only $C_{\varepsilon}^{(1)}>0$ for each parameter $\varepsilon$, but also
the net $C_{\varepsilon}^{(1)}>0$ is monotonically increasing with
$\varepsilon\rightarrow 0$, furthermore (\ref{eqcontained}) holds,
and we are done with this case.
\item $x_2 \notin \widetilde S_{r_1}(x_1)$, that is $\vert x_2-x_1\vert_e<r_1$. Set, for instance, $C_{\varepsilon}^{(1)}=1$.
 For each representative $(x_{\varepsilon}^{(2)})_{\varepsilon}$ of $x_2$ it follows that
 \[
 \vert x_\varepsilon^{(2)}-x_\varepsilon^{(1)}\vert=o(\varepsilon^{\rho_1})
 \]
 and a representative satisfying the desired properties is easily found.
\end{itemize}
Proof of (\ref{lemmatique2}):\\
To show this we consider the asymptotic growth of
$(C_{\varepsilon}^{(1)})_{\varepsilon},\,(C_{\varepsilon}^{(2)})_{\varepsilon},\,\varepsilon^{\rho_1},\,\varepsilon^{\rho_2}$
as well as the monotonicity of $C_{\varepsilon}^{(1)}$. Let $y\in
B_{\leq C_{\varepsilon}^{(2)}\varepsilon^{\rho_2}}(x_{\varepsilon}^{(2)})$.
By the triangle inequality we have that
\begin{equation}\label{est0}
\vert y-x_{\varepsilon}^{(1)}\vert\leq \vert
y-x_{\varepsilon}^{(2)}\vert +\vert
x_{\varepsilon}^{(2)}-x_{\varepsilon}^{(1)}\vert\leq
C_{\varepsilon}^{(2)}\varepsilon^{\rho_2}+\frac
{C_{\varepsilon}^{(1)}\varepsilon^{\rho_1}}{2},
\end{equation}
for all $\varepsilon\in I$. We know further that by the monotonicity $\forall \varepsilon\in I:
C_{\varepsilon}^{(1)}\geq C_{\varepsilon=1}^{(1)}=:C_1$ so that
\begin{equation}\label{est1}
\frac{C_{\varepsilon}^{(2)}}{C_{\varepsilon}^{(1)}}\varepsilon^{\rho_2-\rho_1}\leq
C_1C_{\varepsilon}^{(2)}\varepsilon^{\rho_2-\rho_1}.
\end{equation}
Moreover, since the sharp norm of $C_{\varepsilon}^{(2)}$ equals
$1$, for any $\alpha>0$ we have that
\[
C_{\varepsilon}^{(2)}=o(\varepsilon^{-\alpha}),\;
(\varepsilon\rightarrow 0),
\]
which in conjunction with the fact that $\rho_2>\rho_1$ allows us to
further estimate the right hand side of (\ref{est1}): Obtaining
\[
\frac{C_{\varepsilon}^{(2)}}{C_{\varepsilon}^{(1)}}\varepsilon^{\rho_2-\rho_1}=o(1),\;(\varepsilon\rightarrow
0).
\]
We plug this information into (\ref{est0}). This yields for
sufficiently small $\varepsilon$, say
$\varepsilon<\varepsilon_0^{(1)}$:
\begin{equation}
\vert y-x_{\varepsilon}^{(1)}\vert\leq \frac
{C_{\varepsilon}^{(1)}\varepsilon^{\rho_1}}{2}+ \frac
{C_{\varepsilon}^{(1)}\varepsilon^{\rho_1}}{2}=C_{\varepsilon}^{(1)}\varepsilon^{\rho_1}
\end{equation}
and completes the proof.
\end{proof}
\begin{proposition}\label{propseq}
Any nested sequence of closed balls in $\widetilde{\mathbb K}$
admits a proper sequence of associated euclidean models.
\end{proposition}
\begin{proof}
We proceed step by step so that we can easily read off the inductive argument of the proof in the end.\\
We may assume that for each $i\geq 1$, $r_i>r_{i+1}$. Define $\rho_i:=-\log (r_i)$ (so that $\rho_i<\rho_{i+1}$ for each $i\geq 1$).\\
{\bf Step 1.}\\
Choose a representative $(x^{(1)}_\varepsilon)_\varepsilon$ of $x_1$.\\
{\bf Step 2.}\\
Due to Proposition \ref{lemmatique} (\ref{lemmatique1}) we can
choose a representative $(x^{(2)}_\varepsilon)_\varepsilon$ of $x_2$
and a net  $(C^{(1)}_\varepsilon)_\varepsilon$ of real numbers
satisfying condition (E) such that 
\[
x_\varepsilon^{(2)}\in B_{\leq
\frac{C_\varepsilon^{(1)}\varepsilon^{\rho_1}}{2}}(x_\varepsilon^{(1)})
\]
for all $\varepsilon\in I$.\\
{\bf Step 3.}\\
Similarly, take a representative $(\hat
x^{(3)}_\varepsilon)_\varepsilon$ of $x_3$ and a net  $(\hat
C^{(2)}_\varepsilon)_\varepsilon$ of real numbers satisfying
condition (E) such that such that for each $\varepsilon\in I$
\begin{equation}\label{inclpointball2}
\hat x_\varepsilon^{(3)}\in B_{\leq \frac{\hat
C_\varepsilon^{(2)}\varepsilon^{\rho_2}}{2}}(x_\varepsilon^{(2)}).
\end{equation}
Denote by $\varepsilon_0^{(1)}\in I$ be the maximal $\varepsilon$
such that the inclusion relation
$B^{(2)}_{\varepsilon}\subseteq B^{(1)}_{\varepsilon}$  holds (cf.\ (\ref{lemmatique2}) of Proposition \ref{lemmatique}).
We show now, how to adjust our choice of $\hat x_\varepsilon^{(3)},\;
\hat C_\varepsilon^{(2)}$ such that condition (E) as well as the
inclusion relation (\ref{inclpointball2}) is preserved, however, we
do this in a way such that we moreover achieve the inclusion
relation
\begin{equation}\label{inclrel2balls}
B^{(2)}_{\varepsilon}\subseteq B^{(1)}_{\varepsilon}
\end{equation}
for {\it each} $\varepsilon$. For $\varepsilon < \varepsilon_0^{(1)}$ we leave the choice
unchanged, that is, we set
\[
 x_\varepsilon^{(3)}:=\hat x_\varepsilon^{(3)},\;C_\varepsilon^{(2)}:=\hat C_\varepsilon^{(2)}.
\]
For $\varepsilon \geq \varepsilon_0^{(1)}$, however, we set
\begin{equation}\label{resetstep3}
 x_\varepsilon^{(3)}:=x_\varepsilon^{(2)},\;C_\varepsilon^{(2)}:=\min(\frac{C_\varepsilon^{(1)}}{2}\varepsilon^
 {\rho_1-\rho_2},  \hat C_\varepsilon^{(2)}).
\end{equation}
Therefore, $(C_\varepsilon^{(2)})_\varepsilon$ still satisfies
condition (E), since it is still positive and monotonically
increasing with $\varepsilon\rightarrow 0$. Next, it is evident that
\[
x_\varepsilon^{(3)}\in B_{\leq
\frac{C_\varepsilon^{(2)}\varepsilon^{\rho_2}}{2}}(x_\varepsilon^{(2)}).
\]
still holds for each $\varepsilon\in I$. Finally, by
(\ref{resetstep3}) it follows that the inclusion relation
(\ref{inclrel2balls}) holds now for each $\varepsilon\in I$. For the
inductive proof of the statement one formally proceeds as in Step 3.
Let $k>1$. Assume we have representatives
\[
(x_{\varepsilon}^{(1)})_{\varepsilon},\dots,(x_{\varepsilon}^{(k+1)})_{\varepsilon}
\]
and nets of positive numbers
\[
(C_{\varepsilon}^{(j)})_{\varepsilon}, (1\leq j\leq k),
\]
satisfying condition (E), such that for each $\varepsilon\in I$ we
have:
\[
B_{\leq
C_{\varepsilon}^{(1)}\varepsilon^{\rho_1}}(x_{\varepsilon}^{(1)})\supseteq
B_{\leq
C_{\varepsilon}^{(2)}\varepsilon^{\rho_2}}(x_{\varepsilon}^{(2)})\supseteq\dots\supseteq
B_{\leq
C_{\varepsilon}^{(k-1)}\varepsilon^{\rho_{k-1}}}(x_{\varepsilon}^{(k-1)}),
\]
and for some $\varepsilon_0^{(k-1)}$ we have for each
$\varepsilon<\varepsilon_0^{(k-1)}$
\[
B_{\leq
C_{\varepsilon}^{(k-1)}\varepsilon^{\rho_{k-1}}}(x_{\varepsilon}^{(k-1)})\supseteq
B_{\leq
C_{\varepsilon}^{(k)}\varepsilon^{\rho_{k}}}(x_{\varepsilon}^{(k)}).
\]
Furthermore we suppose the following additional property is
satisfied: For each $\varepsilon\in I$ we have:
\[
x_{\varepsilon}^{(k+1)}\in B_{\leq
\frac{C_{\varepsilon}^{(k)}}{2}\varepsilon^{\rho_k}}(x_{\varepsilon}^{(k)}),
 \]
 where $\rho_k:=-\log r_k$. In the very same manner as above, we can now find a representative $(x_{\varepsilon}^{(k+2)})_{\varepsilon}$ of $x_{k+2}$ and
a net of numbers $(C_{\varepsilon}^{(k+1)})_{\varepsilon}$
satisfying condition (E) such that the above sequential construction
can be enlarged by one ($k\rightarrow k+1$).
\end{proof}
The preceding proposition is a key ingredient in the proof of our
main statement Theorem \ref{theoremsph}:
\begin{proof}
Let $(\widetilde B_i)_{i=1}^{\infty},\; B_i:=\widetilde B_{\leq
r_i}(x_i)\;(i\geq 1)$ be the given nested sequence of dressed balls;
due to Proposition \ref{propseq}, there exists a proper sequence of
associated euclidean models
\[
(B^{(i)}_{\varepsilon})_{\varepsilon}
\]
such that for representatives
$(x^{(i)}_{\varepsilon})_{\varepsilon}$ of $x_i$ ($i\geq 1$) the
above nets are given by
\[
B^{(i)}_{\varepsilon}:=B_{\leq
C_{\varepsilon}^{(i)}\varepsilon^{\rho_i}}(x^{(i)}_{\varepsilon}),\quad\rho_i:=-\log
r_i,\quad C_{\varepsilon}^{(i)}\in\mathbb R_+
\]
for each $(\varepsilon,i)\in I\times\mathbb N$. Since $\mathbb K$ is
locally compact, for each $\varepsilon\in I$ we can choose some
$x_{\varepsilon}\in\mathbb R$ such that
\[
x_{\varepsilon}\in \bigcap_{i=1}^{\infty}B^{(i)}_{\varepsilon}
\]
since for each $\varepsilon\in I$ we have
$B_{\varepsilon}^{(1)}\supseteq
B_{\varepsilon}^{(2)}\supseteq\dots$ . By the construction of the net
$(x_\varepsilon)_\varepsilon$,
we have
\[
\vert x_{\varepsilon}-x_{\varepsilon}^{(i)}\vert \leq
C_{\varepsilon}^{(i)}\varepsilon^{\rho_i}
\]
for each $\varepsilon\in I$. This shows that not only the net $(x_{\varepsilon})_{\varepsilon}$
is moderate (use the triangle inequality), but also gives rise to a
generalized number $x:=(x_{\varepsilon})_{\varepsilon}+\mathcal
N(\mathbb K)$ with the property
\[
\vert x-x_i\vert_{\rm e}\leq r_i
\]
for each $i$. Therefore we have that
\[
x\in\bigcap_{i=1}^{\infty}\widetilde B_i\neq \emptyset
\]
which yields the claim: $\widetilde{\mathbb K}$ is spherically
complete.
\end{proof}
\section{A Hahn-Banach Theorem} Let $L$ be a subfield of
$\widetilde{\mathbb K}$. Let $(E,\|\cdot\|)$ be an ultrametric normed $L$-linear space. We
call $\varphi$ an $L$- linear functional on $E$, if $\varphi$ is an
$L$- linear mapping on $E$ with values in $\widetilde{\mathbb K}$.
$\varphi$ is continuous if and only if
\[
\|\varphi\|:=\sup_{0\neq x\in E}\frac{\vert
\varphi(x)\vert_{\rm e}}{\|x\|}<\infty.
\]
We denote the space of all continuous $L$-linear functionals on $E$ by $E'_L$.
\begin{remark}\rm
Note that nontrivial subfields $L$ of $\widetilde{\mathbb K}$ exist.
For instance, one can choose $\mathbb K(\alpha)$ with
$\alpha=[(\varepsilon)_{\varepsilon}]\in\widetilde{\mathbb K}$ or
its completion with respect to $\vert \;\;\vert_{\rm e}$, the Laurent series
over $\widetilde{\mathbb K}$. Moreover, given an ultra-pseudo-normed $\widetilde{\mathbb C}$-
module $(\mathcal G,\mathcal P)$, the $L$-linear space $E$ generated by elements of $\mathcal G$
is an an ultrametric normed $L$-linear space.
\end{remark}
Having introduced these notions we show that the following version of
the Hahn-Banach Theorem holds:
\begin{theorem}
Let $V$ be an $L$-linear subspace of $E$ and $\varphi\in V'_L$. Then
$\varphi$ can be extended to some $\psi\in E'_L$ such that
$\|\psi\|=\|\varphi\|$.
\end{theorem}
\begin{proof}
We follow the lines of the proof of Ingleton's theorem (cf.\
\cite{Ingleton}) in the fashion of (\cite{zAR}, pp.\ 194--195). To
start with, let $V$ be a strict $L$-linear subspace of $E$ and let
$a\in E\setminus V$. We first show that $\varphi\in V'_L$ can be
extended to $\psi\in (V+La)'_L$ under conservation of its norm. To
do this it is sufficient to prove that such $\psi$ satisfies for
each $x\in V$:
\begin{eqnarray}\label{ineqnorm}
\|\psi(x-a)\|&\leq&\|\psi\|\cdot\|x-a\| \\\nonumber
\|\varphi(x)-\psi(a)\|&\leq&\|\varphi\|\cdot\|x-a\|=:r_x.
\end{eqnarray}
To this end define for each $x$ in $V$ the dressed ball
\[
B_x:=B_{\leq r_x}(\varphi(x)).
\]
Next we claim that the family $\{B_x\mid x\in V\}$ of dressed balls
is nested. To see this, let $x,y \in V$. By the linearity of
$\varphi$ and the ultrametric (strong) triangle inequality we have
\[
\vert \varphi(x)-\varphi(y)\vert_{\rm e}\leq \|
\varphi\|\cdot\|x-y\|\leq\|\varphi\|\max(\|x-a\|,\|y-a\|)=\max(r_x,r_y).
\]
Therefore we have $B_x\subseteq B_y$ or $B_y\subseteq B_x$. According to Theorem \ref{theoremsph}, $\widetilde{\mathbb
K}$ is spherically complete, therefore we can choose
\[
\alpha\in \bigcap_{x\in V}B_x
\]
and further define $\psi(a):=\alpha$. According to (\ref{ineqnorm}) we therefore have for each $z\in V$ and for each $\lambda \in L$ \footnote{Note that since $E$ is a normed $L$-linear space,
the restriction of $\vert \cdot\vert_{\rm e}$ to $L$ is multiplicative.},
\[
\vert\psi(z-\lambda a)\vert_{\rm e}=\vert\lambda\vert_{\rm e}\cdot\vert
\psi(z/\lambda-a)\vert_{\rm e}\leq \vert \lambda\vert_{\rm e}\, r_{z/\lambda}=\vert
\lambda\vert_{\rm e}\cdot \|\varphi\|\cdot\|z/\lambda-a\|=\|\varphi\|\cdot
\|z-\lambda a\|
\]
which shows that $\psi$ is an extension of $\varphi$ onto $V+La$ and
$\|\psi\|=\|\varphi\|$.

The rest of the proof is the standard one-an application of Zorn's
Lemma.
\end{proof}
Let $(\mathcal G, \|\cdot\|)$ be an ultra-pseudo-normed $\widetilde{\mathbb K}$ module and denote by
$\mathcal L(\mathcal G,\widetilde{\mathbb K})$ the space of continuous linear functionals on $E$ (according to the notation in \cite{Garetto2,Garetto1}). We end this section
by posing the following conjecture:
\begin{conjecture}
Let $\mathcal V$ be a submodule of $\mathcal G$ and let $\varphi\in \mathcal L(\mathcal V,\widetilde{\mathbb K})$. Then
$\varphi$ can be extended to some element $\psi\in\mathcal L(\mathcal G,\widetilde{\mathbb K})$ such that
$\|\psi\|=\|\varphi\|$.
\end{conjecture}
\section*{Appendix} Finally, it is worth
mentioning that apart from the standard Fixed Point Theorem due to
Banach, a non-archimedean version is available in spherically
complete ultrametric spaces (therefore, also on $\widetilde{\mathbb
K}$, cf.\ \cite{priess1}, and for a recent generalization cf.\
\cite{priess2}):
\begin{theorem}
Let ($M,d$) be a spherically complete ultrametric space and $f:
M\rightarrow M$ be a mapping having the property
\[
\forall x,y\in M: d(f(x),f(y))<d(x,y).
\]
Then $f$ has a unique fixed point in $M$.
\end{theorem}
\section*{Acknowledgment}
I am indebted to Prof.\ M.\ Kunzinger (Vienna) and Prof.\ S.\
Pilipovi{\'c} (Novi Sad) for reading carefully the manuscript and for
helpful advise. Also, I thank the anonymous referee for valuable and detailed suggestions which lead
to the final form of this paper.


\begin{thebibliography}{10}

\bibitem{DHPV}
{\sc A.~Delcroix, M.~F. Hasler, S.~Pilipovi{\'c}, and V.~Valmorin}, {\em
  Generalized function algebras as sequence space algebras}, Proc. Amer. Math.
  Soc., 132 (2004), pp.~2031--2038 (electronic).

\bibitem{Scarpalezos1}
{\sc A.~Delcroix and D.~Scarpalezos}, {\em Sharp topologies on ($\widetilde
  {C}$,$\widetilde {E}$,$\widetilde {P}$)-algebras}, in Nonlinear theory of
  generalized functions (Vienna, 1997), vol.~401 of Chapman \& Hall/CRC Res.
  Notes Math., Chapman \& Hall/CRC, Boca Raton, FL, 1999, pp.~165--173.

\bibitem{Garetto2}
{\sc C.~Garetto}, {\em Topological structures in {C}olombeau algebras:
  investigation of the duals of {$\mathcal G_c(\Omega),\, \mathcal G(\Omega)$
  and $\mathcal G_{\mathcal S}(\mathbb R^n)$}}, Monatsh. Math., 146 (2005),
  pp.~203--226.

\bibitem{Garetto1}
\leavevmode\vrule height 2pt depth -1.6pt width 23pt, {\em Topological
  structures in {C}olombeau algebras: topological {$\widetilde{\mathbb
  C}$}-modules and duality theory}, Acta Appl. Math., 88 (2005), pp.~81--123.

\bibitem{Bible}
{\sc M.~Grosser, M.~Kunzinger, M.~Oberguggenberger, and R.~Steinbauer}, {\em
  Geometric theory of generalized functions with applications to general
  relativity}, vol.~537 of Mathematics and its Applications, Kluwer Academic
  Publishers, Dordrecht, 2001.

\bibitem{Ingleton}
{\sc W.~Ingleton}, {\em The {H}ahn-{B}anach theorem for non-archimedean valued
  fields}, Proc. Cambridge Phil. Soc., 48 (1952), pp.~41--45.

\bibitem{Eyb4}
{\sc E.~Mayerhofer}, {\em On the characterization of p-adic
  {C}olombeau-{E}gorov generalized functions by their point values}, to appear
  in Math. Nachr.,  (2006).

\bibitem{MO1}
{\sc M.~Oberguggenberger and M.~Kunzinger}, {\em Characterization of
  {C}olombeau generalized functions by their pointvalues}, Math. Nachr., 203
  (1999), pp.~147--157.

\bibitem{priess1}
{\sc S.~Prie{\ss}-Crampe}, {\em Der {B}anachsche {F}ixpunktsatz f\"ur
  ultrametrische {R}\"aume}, Results Math., 18 (1990), pp.~178--186.

\bibitem{priess2}
{\sc S.~Priess-Crampe and P.~Ribenboim}, {\em Fixed point and attractor
  theorems for ultrametric spaces}, Forum Math., 12 (2000), pp.~53--64.

\bibitem{zAR}
{\sc A.~M. Robert}, {\em A course in {$p$}-adic analysis}, vol.~198 of Graduate
  Texts in Mathematics, Springer-Verlag, New York, 2000.

\end{thebibliography}
\end{document}